\documentclass[12pt]{amsart}
\usepackage[english]{babel}

\newcommand{\Z}{\mathbb Z}

\title{Comments on Sen's restricted preferences and rational choice} 

\author{Labib Haddad}
\address{120 rue de Charonne, 75011 Paris, France}
\email{labib.haddad@wanadoo.fr}

\usepackage{amssymb}
\usepackage{amsmath}

\newcommand{\su}{\subsection*}
\newcommand{\head}{\section*}
\newcommand{\noi}{\noindent}

\newcommand{\cal}{\mathcal}

\newcommand{\stm}{\smallsetminus}

\begin{document}
\maketitle

\thispagestyle{empty}

\markboth{Labib Haddad}{Condorcet paradox}

\hfill{\it Out, out, brief candle!}

\hfill{Macbeth, Act V, Scene 5}

\

By now, most mathematicians have  probably heard about voting paradoxes. One of the most simple, and easy to describe, among those paradoxes  is the following.

In certain situations, an assembly can very well decide, by successive majority votes,  that $a$ is preferable to $b$, that $b$ is preferable to $c$, but that $c$ is preferable to $a$. No matter what $a$, $b$, and $c$, are, this is a cyclic, paradoxical set of decisions: Those $a, b, c$, might be candidates in some election, or options, choices to be compared! Or what else you wish and can figure out!

For brevity's sake, one writes $a > b > c > a$.

This paradox has been noticed a long time ago; one way to speak about it, in French, is to say \lq\lq l'effet Condorcet".

In most decision procedures, this situation is inevitable. Attempts have been made to find a way out and to avoid it, with no avail!

Many details can be found in the book  [4] by Sen.

Before we go on, just a few words about Guilbaud's paper, [2], concerning the problem of aggregation of decisions and voting systems.

\

\head{On Guilbaud's voting systems}

\

Studying  Condorcet's Paradox, Guilbaud, [2], devised a schema for a whole set of possible voting systems.

Let $A$ be an assembly. Each subset $K$ of $A$ is called {\bf a coalition}. Let the complement, $K^c = A \stm K$, denote the opposing coalition to $K$.

\

A {\bf voting system} for the assembly $A$ is a  set $\cal E$ of coalitions of $A$ (to be used as substitutes for ordinary majorities). In other words, $\cal E$ is a subset of the set $S(A)$ of all subsets of the set $A$. The elements of $\cal E$ are, by definition, the efficient coalitions.
There are two requirements put on $\cal E$.

\

{\bf C1} A coalition is efficient if and only if the opposing coalition is not.

\

{\bf C2} Any coalition that contains an efficient coalition is efficient, itself.

\

The system ($A,\cal E$) operates in the following way. Questions submitted to the assembly are dichotomous choices which the members of the assembly answer by {\bf yes} or {\bf no}. Thus two coalitions form, the {\bf pros} and the {\bf cons}. The two coalitions are complementary. 
Of these two coalitions the efficient one, the one which belongs to $\cal E$, outweighs the other: the decision is hers!

Of course, the majority system is one of Guilbaud's voting systems provided a {\it special clause} is added to it, such as a cast vote,  in case there is an even number of voters.
By the way, dictatorship is also a Guilbaud voting system!

\

Now, to make a long story short, there is a problem which stated in terms of pure mathematics is as follows.

\

\head{The problem}

\

Suppose there are three candidates, $a,b,c$, to be compared, two by two,  by successive votes of an assembly $A$  according to a Guilbaud voting system, $(A,\cal E)$. Each member of the assembly having his  own preferences, to be called his {\bf profile}, one among the six following linear orders:

\

$1 \ \ a>b>c$ 

$2 \ \ a>c>b$ 

$3 \ \ c>a>b$ 

$4 \ \ c>b>a$ 

$5\ \  b>c>a$ 

$6 \ \ b>a>c$.

\

\noi Those profiles are labeled by the elements $\{1,2,3,4,5,6\}$ of the cyclic group $\Z/6\Z$ of order $6$.

The labelling has the following peculiarities. Rankings $p$ and $p + 1$ always have either a first same candidate or a last same candidate. To go from $p$ to $p + 1$, one \lq\lq disturbs" as little as possible these rankings, i.e., one simply swaps the ranks of either the two first or of the two last candidates. In a sense, rankings $p$ and $p + 1$ are \lq\lq as close as possible" to each other. Finally, rankings $p$ and $p + 3$ are \lq\lq opposite". For example, one has

\

$1 \ \ a>b>c$ 

$4 \ \ a<b<c$.

\

\noi Thus, in particular, the same candidate occupies the second rank in two different rankings, $p$ and $q$, if and only if we have $p + q = 3$, i.e., when the two rankings $p$ and $q$ are opposite.

\

Let $K(p)$ be the set of the members of the assembly whose profile is $p$. Similarly, denote $K(p,q)$ the union of coalitions $K(p$) and $K(q)$, and let $K(p,q,r)$ be the union of coalitions $K(p)$, $K(q)$, and $K(r)$.

\su{Problem} Find conditions on those coalitions such that the final decision of the assembly be one of the six rankings, not a cyclic ranking!

\

\su{Answer} Condition (C) that there exists a profile $p$ such that  both coalitions $K(p, p + 1, p + 2)$ \ and $K(p+1,p+2,p+3)$ are efficient,
is necessary and sufficient for the decision of the assembly to be one of the six rankings, not cyclic.

\

\noi {\bf Indeed}, suppose that  $K(p, p + 1, p + 2)$ and $K(p+1,p+2,p+3)$ are both efficient. The ranking $p$ has the form $x > y > z$.

\

There are only two possible cases.

\

\noi In the first case, one has

$p \quad \quad \ x>y>z$ 

$p+1 \ \  x>z>y$ 

$p+2 \ \ z>x>y$ 

$p+3 \ \ z>y>x$.

\noi The efficient coalition $K(p, p + 1, p + 2)$ thus imposes the collective preference $x > y$ and efficient coalition $K(p+1, p+2, p+3)$ the collective preference $z > y$. Whatever the collective preference between $x$ and $z$, the collective ranking will always be linear! [It will be either $p+1$ or $p+2$.]

\noi In the second case, one has

$p \quad \quad \ x>y>z$ 

$p+1 \ \ y>x>z$ 

$p+2 \ \ y>z>x$ 

$p+3 \ \ z>y>x$.

 \
 
\noi Here, the two efficient coalitions $K(p, p + 1, p + 2)$ and $K(p + 1,p+2,p+3)$ force, respectively, collective preferences $y > z$ and $y > x$ so that the collective ranking is  linear. [This will, again, be either p+1 or p+2.] 

\

\noi Now, conversely, suppose the collective ranking of the three candidates is linear. Even if we have to change the names of the candidates, we can assume that this ranking is $a > b > c$. The coalition of voters for whom $a > b$ is efficient: This coalition is none other than $K(1, 2, 3)$. Similarly, the coalition of voters for whom $a > c$ is efficient and that is coalition $K (6, 1, 2)!$ So that condition (C) is satisfied.  

Condition (C)  is both necessary and sufficient! \qed

\noi Of course this applies to the majority system (provided a special clause is added to it). 

\head{A quotation}

\

Here is a long quotation from the paper, [1], by Elsholtz and List.

\

\lq\lq Condorcet's paradox shows that pairwise majority voting over three or more candidates can lead to cyclical majority preferences, even when the preferences of individual voters are transitive ...

... A large literature in social choice theory addresses the threat posed by cycles for the functioning of democratic decision mechanisms ... Several sufficient conditions for the avoidance of cycles have been identified. Black   [identified one such] condition ... single-peakedness. Later, other sufficient conditions for transitive majority preference were found, amongst them single-cavedness, separability into two groups, and latin-squarelessness [for instance].

In a famous paper, Sen generalized these results, showing that a condition that is less demanding than, but implied by, each of these conditions is already sufficient for avoiding cycles. Sen's condition is called {\it triplewise value-restriction}. However, Sen's condition and theorem are not intuitively straightforward. This note aims to make the mechanism underlying Sen's result easily accessible, by giving an elementary proof of Sen's theorem, together with a simple reformulation of the condition of triplewise value-restriction.... Finally, we suggest that, although there is still some logical space between Sen's sufficient condition for the avoidance of cycles and a necessary and sufficient condition, this space may be too narrow to allow an appealing generalization of Sen's condition."

It might be interesting to compare that quotation with the above elaboration! We must say that we have already wrote similar comments about the paper by Elsholtz and List, in [3].

\

\

\head{References}

\

\noi[1]{\sc ELSHOLTZ Christian and LIST Christian}, {\sl A Simple Proof of Sen's Possibility Theorem on Majority Decisions}, Elemente der Mathematik, {\bf 60} (2005) no 2, 45-56.

\

\noi [2] {\sc GUILBAUD G.Th.}, {\sl Les th\Žories de l'int\Žr\t g\Žn\Žral et le probl\me logique de l'agr\Žgation}, \' Economie appliqu\Že, {\bf 5} (1952) no 4, oct.-d\Žc., 501-551. 

\

\noi[3] {HADDAD Labib}, {\sl Un outil incomparable \ : l'ultrafiltre}.
Tatra Mt. Math. Publ., {\bf 31}  (2005) 131-176. See also a translation in English,  {\sl The ultrafilter: A peerless tool},
{\tt arXiv:math/0702587v2, 5 May 2021}.

\

\noi [4] {\sc SEN Amartya K.}, {\sl Choice, welfare, and measurement}, (1997) Cambridge, MA: Harvard University Press 

\

\enddocument